\documentclass[11pt]{article}    
\usepackage{amstext}
\usepackage{amsmath}
\usepackage[dvips]{graphicx}
\usepackage{amsfonts}
\usepackage{amssymb}
\usepackage{amsxtra}
\usepackage{natbib}

\newtheorem{lemma}{Lemma}[section]
\newtheorem{remark}{Remark}[section]
\newtheorem{cor}{Corollary}[section]
\newtheorem{proposition}{Proposition}[section]
\newenvironment{proof}{\textsc{Proof:}}{\mbox{ } \hfill $\Box$ \vspace{2mm}}

\numberwithin{equation}{section}

\newcommand{\Fb}{\mathbb{F}}
\newcommand{\R}{\mathbb{R}}
\newcommand{\N}{\mathbb{N}}

\renewcommand{\P}{\mathbb{P}}
\newcommand{\E}{\mathbb{E}}

\newcommand{\pk}{\Phi^{(k)}}

\newcommand{\eps}{\varepsilon}
\newcommand{\p}{\widetilde{\Phi}}
\newcommand{\px}{\widetilde{\Phi}^{(0)}}
\newcommand{\py}{\widetilde{\Phi}^{(1)}}
\newcommand{\Ei}{\mathbb{E}^{\phi_0,\phi_1}}
\newcommand{\Co}{\mathbb{C}_0}
\newcommand{\bGamma}{\mathbf{\Gamma}}
\newcommand{\C}{\mathbf{C}}

\renewcommand{\epsilon}{\varepsilon}

\newcommand{\wPhiz}{\widetilde{\Phi}^{(0)}}
\newcommand{\wPhio}{\widetilde{\Phi}^{(1)}}

\newcommand{\wbPhi}{\widetilde{{\Phi}}}

\title{\LARGE \bf
\sc  Optimal Time to Change Premiums
\thanks{This work was partially supported by the US Army Pantheon Project and National Science Foundation under grant DMS-0604491}
\thanks{The authors would like to thank the participants of the
Statistics and Operations Research Seminar at the University of
North Carolina at Chapel Hill, and the Industrial Engineering and
Operations Research Seminar at the University of California,
Berkeley.}}

\author{Erhan Bayraktar 
\thanks{E. Bayraktar is with the Department of Mathematics, University
of Michigan, Ann Arbor, MI 48109, USA, email: erhan@umich.edu}
 \and  H. Vincent Poor%
\thanks{H. V. Poor is with the School of Engineering and Applied Science, Princeton University,
        Princeton, NJ 08544, USA, email: poor@princeton.edu}%
}
\date{}

\begin{document}
\maketitle

\begin{abstract}
The claim arrival process to an insurance company is modeled by a
compound Poisson process whose intensity and/or jump size
distribution changes at an unobservable time with a known
distribution.  It is in the insurance company's interest to detect
the change time as soon as possible in order to re-evaluate a new
fair value for premiums to keep its profit level the same.  This
is equivalent to a problem in which the intensity and the jump
size change at the same time but the intensity changes to a random
variable with a know distribution. This problem becomes an optimal
stopping problem for a Markovian sufficient statistic. Here, a
special case of this problem is solved, in which the rate of the
arrivals moves up to one of two possible values, and the Markovian
sufficient statistic is two-dimensional.
\end{abstract}

\section{Introduction}

\maketitle

In insurance risk theory, the claim arrivals are modeled by a
compound Poisson process. The total claim up to time $t$ is given
by
\begin{equation}\label{eq:claim}
X_t=X_0+ \sum_{k=1}^{N_t}Y_k, \quad t \geq 0,
\end{equation}
where the number of claims up to time $t$, $N_t$, is a Poisson
process with intensity $\lambda_0$. The claim size process
$(Y_k)_{k \in \mathbb{N}}$ is assumed to consist of independent
and identically distributed $R^{d}$ valued random variables with
distribution function $\nu_0$. In order to compensate for the
liabilities the insurance company has to pay out, it collects
premiums at a such rate that it has a fair chance of survival.

In this paper, we will study the model in (\ref{eq:claim}) with
two types of regime shift. At time $\theta^a$ the intensity of the
Poisson process changes from $\lambda_0$ to $\lambda_1$, and at
time $\theta^b$, the distribution of the claim size changes from
$\nu_0$ to $\nu_1$. (These measures are assumed to be absolutely
continuous with respect to each other.) Both $\theta^a$ and
$\theta^b$ are unknown at time 0, and they are unobservable. It is
in the insurance company's interest to detect \emph{the change
time or the disorder time} $\theta \triangleq \theta^a \wedge
\theta^b =\min\{\theta^a,\theta^b\}$ as soon as possible and to
re-evaluate a new fair value for premiums in order to keep the
profit level the same.

We assume that the times of regime shift are independent of each
other and that they have an exponential prior distribution
\[
 \quad \P\{\theta^i>t\}=
(1-\pi^i)e^{-\lambda^i t}, \quad i \in \{a,b\}, \quad t\geq 0,
\]
for $\lambda^i>0$. At time $\theta$, we do not know what the
intensity is for sure: it is either $\lambda_0$ (a change has
occurred in the distribution of the claim size) or $\lambda_1$ (a
change occurred in the intensity). In fact at time $\theta$, the
value of intensity changes from $\lambda_0$ to the random variable
$\Lambda$ where
\begin{equation}\label{eq:Lambda}
\Lambda =
\begin{cases}
\lambda_1 & \text{with probability} \quad
\frac{\lambda^a}{\lambda^a+\lambda^b}
\\ \lambda_0 & \text{with probability} \quad\frac{\lambda^b}{\lambda^a+\lambda^b}.
\end{cases}
\end{equation}

At time $\theta$ the distribution of the claim size changes from
$\nu_0$ to $\nu$, where
\begin{equation}\label{eq:nu-convex}
\nu=\frac{\lambda^a}{\lambda^a+\lambda^b} \nu_0+
\frac{\lambda^b}{\lambda^a+\lambda^b} \nu_1.
\end{equation}
Now consider a related more general problem in which at the
disorder time $\theta$ the compound process introduced in
(\ref{eq:claim}) changes its intensity from $\mu \in \mathbb{R}_+$
to a random variable $\Lambda$ (at first we will first allow the
distribution of this random variable to be as general as possible)
and the distribution of the claim sizes change from $\beta_0$ to
$\beta_1$ (these two measures are assumed to be absolutely
continuous with respect to each other). The distribution of
$\theta$ is given by
\begin{equation}\label{eq:dist-theta}
\P\{\theta=0\}=\pi, \quad \P\{\theta>t| \theta>0\}=e^{-\lambda t},
\,t \geq 0.
\end{equation}
The random variables $\Lambda$ and $\theta$ are independent.

In this more general problem the aim is to detect the unknown and
unobservable time $\theta$ as quickly as possible given the
observations from the incoming claims. More precisely, we would
like to find a stopping time $\tau$ of the observation process
that minimizes the penalty function
\begin{equation}\label{eq:penalty}
R_{\tau}(\pi) \triangleq \P\{\tau<\theta\}+ c\,\E[\tau-\theta]^+,
\end{equation}
which is the sum of the frequency of $\mathbb{P}(\tau<\theta)$
false alarms and the expected cost $c\,
\mathbb{E}\left[(\tau-\theta)^{+}\right]$ of detection delay.

We are interested in solving this more general problem for three
reasons. First, setting $\pi=0$, $\lambda=\lambda^a+\lambda^b$,
$\mu=\lambda_0$, $\beta_0=\nu_0$ and $\beta_1=\nu$, and the
distribution of $\Lambda$ to be the Bernoulli distribution in
(\ref{eq:Lambda}) we see that solving this more general problem
also leads to a solution of the main problem introduced in the
second paragraph. Second, in the general problem if we set
$\Lambda$ to be a constant, then we obtain a version of the main
problem in which the rate change and change of the distribution of
the claim sizes occur simultaneously. This case was analyzed by
\cite{ds} and \cite{gapeev}.  Finally, the more general problem
represents a situation in which the insurance company has only
some apriori information about the post disorder rate $\lambda_1$,
but the company can not pin $\lambda_1$ down to a constant because
it might only have very few claims after the regime change occurs.
In fact, the company wants to detect the regime change as soon as
possible, so there is not really any time to collect data to
estimate $\lambda_1$. This change detection problem when the
underlying process $X$ is a (simple) Poisson process was recently
analyzed by \cite{bdk05}. This corresponds to setting
$\beta_0=\nu_0$ and $\beta_1=\nu_0$ in the current setting.

The compound/simple Poisson disorder problem is one of the rare
instances in which a stochastic control problem with partial
information can be handled. The (simple) Poisson disorder problem
with linear penalty for delay was partially solved by
\cite{galchuk71}, \cite{davis76} and \cite{wdavis}. This problem
later was solved by \cite{MR2003i:60071}. \cite{BD03} solved the
simple Poisson disorder problem for exponential penalty for delay,
and \cite{BD04} solved the standard Poisson disorder problem.
These results were recently extended by \cite{ds} (using the
results developed in \cite{bdk05}) and \cite{gapeev} for compound
Poisson procesesses. On the other hand \cite{bdk05} solved the
simple Poisson disorder problem when the post disorder rate is a
random variable and \cite{bs} solved this problem for the case
with a Phase-type disorder distribution.

We will first show that our problem is equivalent to an optimal
stopping problem for a Markovian sufficient statistic. As in
\cite{bdk05} it turns out that the dimension of the sufficient
statistic is finite dimensional if the distribution of the random
variable $\Lambda$ is discrete
 with finitely many atoms. We will study the case of
a binary distribution in more detail. In particular, we will
analyze the case when the post-disorder rate only goes up. We are
able to show that the intuition that a decision would sound the
alarm only at the times when it observes an arrival does not in
general hold, see Remark~\ref{rem:only-at-jump-times}. This
intuition becomes relevant only when $\lambda$ and $c$ are small
enough, i.e. when the disorder intensity and delay penalty are
small. By performing a sample path analysis we are able to find
the optimal stopping time exactly for most of the range of
parameters. For the rest of the parameter range we provide upper
and lower bounds on the optimal stopping time. To show the
existence of the optimal stopping problem for the cases when we
can not determine it exactly we make use of the characterization
of the value function of the optimal stopping time as the fixed
point of a functional operator, as in \cite{bdk05}. We use this
approach since the free boundary problems associated with our
problem turns out to be quite difficult to manage as it involves
integro-differential equations and the failure of the smooth fit
principle is expected. This characterization can be used to
calculate the value function through an iterative procedure. From
this characterization we are able to infer that the free
boundaries are decreasing convex curves located at the corner of
$\R_+^2$. Using our sample path analysis, we are able to determine
a certain subset of the free boundary exactly.

The rest of the paper is organized as follows. In
Section~\ref{sec:reference-prob}, we give a more precise
probabilistic description of the disorder problem and introduce a
reference probability measure $\P_0$ under which the observations
are coming from a compound Poisson process whose jump distribution
does not change over time. In Section~\ref{sec:mss}, we show that
the disorder problem can be transformed into an optimal stopping
problem for a Markovian sufficient statistic. The Markovian
sufficient statistic may not be finite dimensional and we show in
this section that it is finite dimensional when the distribution
of the post disorder rate has finitely many atoms. In
Section~\ref{sec:Bernoulli}, we  find the autonomous sufficient
statistic for any Bernoulli distribution. Also we set up an
optimal stopping problem for a Bernoulli sufficient statistic
 when the post disorder rate can only move up. Section~\ref{sec:bounds} contains some of our main results in
 which by performing a sample path analysis we either find the optimal stopping time exactly or provide upper and lower bounds.
 We also show that
the optimal stopping time is finite $\mathbb{P}_0$-almost surely.
Section~\ref{sec:optimal-stopping-time} provides a useful
characterization of the value function as a limit of a sequence of
other value functions. Since the proofs of the results in this
section are similar to the ones in \cite{bdk05} we omit them,
except the result in which we show that the optimal stopping time
we constructed is the smallest optimal stopping time and a few
other that we prefer to keep for readers convenience.

\section{A Reference Probability Measure}\label{sec:reference-prob}

We will first introduce a reference probability measure $\P_0$
under which the observations have a simpler form, namely they come
from a compound Poisson process whose rate and jump distribution
do not change over time. Next, we will construct the model that we
briefly described in the introduction in the paragraph before
(\ref{eq:dist-theta}).

Let us start with a probability space $(\Omega, \mathcal{F},
\P_0)$ and consider a standard Poisson process $N=\{N_t: t \geq
0\}$ with rate $\mu$; independent and identically distributed
strictly positive random variables $Y_1, Y_2,...$ with a common
distribution $\beta_0$ on $\mathbb{R}^d$ independent of the
Poisson process; a random variable $\theta$ independent of the
previously described stochastic elements on this probability space
whose distribution is given by
\begin{equation}
\P_0\{\theta=0\}=\pi, \quad \P_0\{\theta>t| \theta>0\}=e^{-\lambda
t}, \,t \geq 0;
\end{equation}
a random variable $\Lambda$ independent of the other stochastic
elements whose distribution is $\gamma(\cdot)$. This distribution
charges only the positive real numbers. We will assume that
\begin{equation}\label{eq:moment}
m^{(k)} \triangleq \int_{\mathbb{R}_{+}} (v-\mu)^{k}
\gamma(dv)<\infty, \quad k \in \mathbb{N}_0.
\end{equation}
 Let the process $X=\{X_t: t \geq 0\}$
be the compound Poisson process defined as in (\ref{eq:claim}) and
$\mathbb{F}=\{\mathcal{F}_t\}_{t \geq 0}$ be the natural
filtration of $X$. We will also define an initial enlargement of
$\mathbb{F}$, $\mathbb{G}=\{\mathcal{G}_t\}_{t \geq 0}$ by setting
$\mathcal{G}_t \triangleq \mathcal{F}_t \vee \sigma\{\theta,
\Lambda\}$. $\mathcal{G}_t$ is the information available to a
\emph{genie} at time $t$ that also observes the realizations of
the disorder time $\theta$ and post-disorder rate $\Lambda$. Let
$\beta_1 (\cdot)$ be a probability measure on $\mathbb{R}^d$ which
is absolutely continuous with respect to $\beta_0 (\cdot)$. We
will denote by $r$ the Radon-Nikodym derivative
\begin{equation}
r(y) \triangleq \frac{d \beta_1}{d \beta_0} (y), \quad y \in
\mathbb{R}^d.
\end{equation}
The process
\begin{equation}\label{eq:Z}
Z_t \triangleq \frac{L_t}{L_{\theta}}1_{\{\tau \leq t\}} +
1_{\{\tau>t\}}, \quad t \geq 0,
\end{equation}
is a $\mathbb{G}$-martingale where
\begin{equation}
L_t \triangleq e^{-(\Lambda-\mu)t}
\prod_{k=1}^{N_t}\left[\frac{\Lambda}{\mu} r(Y_k)\right].
\end{equation}
The positive martingale $Z$ defines a new probability measure
$\mathbb{P}$ on every $(\Omega, \mathcal{G}_t)$, $t \geq 0$ by
\begin{equation}
\frac{d \mathbb{P}}{d \mathbb{P}_0}\bigg|_{\mathcal{G}_t}= Z_t,
\quad t \geq 0.
\end{equation}
Note that since $Z_0=1$, $\mathbb{P}$ and $\mathbb{P}_0$ agree on
$\mathbb{G}_0=\sigma\{\theta, \Lambda\}$, i.e. the random
variables $\theta$ and $\Lambda$ are independent and have the same
distribution under both $\mathbb{P}$ and $\mathbb{P}_0$. On the
other hand using the Girsanov Theorem for jump processes (see e.g
\cite{cont}, \cite{ds}) we conclude that the process $X$ is a
$(\mathbb{P}, \mathbb{G})$-compound Poisson process whose arrival
rate $\mu$ and jump distribution $\beta_0$ changes at time
$\theta$ to $\Lambda$ and $\beta_1$, respectively. In other words,
on the probability space $(\Omega, \mathcal{F}, \mathbb{P})$,  we
have exactly the model posited in the Introduction section in the
paragraph between (\ref{eq:nu-convex}) and (\ref{eq:dist-theta}).

\section{Markovian Sufficient Statistics}\label{sec:mss}

In this section, we will show that the stopping problem posed in
(\ref{eq:penalty}) can be formulated as an optimal stopping
problem for a Markovian sufficient statistic, which is in general
infinite dimensional. In the following sections we will see that
depending on the structure of the prior of $\Lambda$ the
sufficient statistic can be finite dimensional.

Let us denote all the $\mathbb{F}$-stopping times by $\mathcal{S}$
and introduce the $\mathbb{F}$-adapted processes
\begin{equation}\label{eq:defn-phik}
\Pi_t \triangleq \mathbb{P}\{\theta \leq t| \mathcal{F}_t\}, \quad
\text{and} \quad \Phi_t^{(k)} \triangleq
\frac{\E\left[(\Lambda-\mu)^{k}1_{\{\theta \leq
t\}}|\mathcal{F}_t\right]}{1-\Pi_t}, \quad k\in \mathbb{N}, t\geq
0.
\end{equation}
$\Pi_t$ is the a posteriori probability process and is the updated
probability that the disorder happened at or before time $t$ given
all the information up to time $t$. $\pk$ can be read as an
\emph{odds-ratio process}, and in fact
$\Phi^{(0)}=\frac{\Pi_t}{1-\Pi_t}$.

Using Proposition 2.1 in \cite{BD04} we can write the Bayes error
in (\ref{eq:penalty}) as
\begin{equation}\label{eq:new-penalty}
R_{\tau}(\pi)=1-\pi+
c(1-\pi)\mathbb{E}_0\left[\int_0^{\tau}e^{-\lambda
t}\left(\Phi_t^{(0)}-\frac{\lambda}{c}\right)dt\right], \quad \tau
\in \mathcal{S},
\end{equation}
where the expectation $\mathbb{E}_0$ is taken under the reference
probability measure $\mathbb{P}_0$. As we can see from
(\ref{eq:new-penalty}), finding an optimal stopping time for the
quickest detection problem would be considerably easier if the
process $\Phi^{(0)}$ is Markovian and its natural filtration
coincides with the filtration generated by the observations. In
that case we would just have to solve a one-dimensional optimal
stopping problem. This is not true, however, unless $\Lambda$ has
only one possible value to take. The following lemma shows that
the whole sequence $\{\pk\}_{k \in \mathbb{N}}$ is a Markovian
sufficient statistic for our detection problem. This result also
will help us develop sufficient conditions under which a finite
dimensional sufficient statistic exists.

\begin{lemma}\label{lem:dyn-of-Phi}
Let $m^{(k)}$ be as in (\ref{eq:moment}). Then the dynamics of
$\Phi^{(k)}$ can be written as
\begin{equation}\label{eq:dyn-Phi}
d\Phi^{(k)}_t=(\lambda(m^{(k)}+\Phi_t^{(k)})-\Phi_t^{(k+1)})dt+\Phi_{t-}^{(k)}\int_{y
\in \mathbb{R}^{d}}(r(y)-1)p(dt dy)+
\Phi_{t-}^{(k+1)}\frac{1}{\mu}\int_{y \in \mathbb{R}^{d}}r(y)p(dt
dy),
\end{equation}
with $\Phi^{(k)}_0=\frac{\pi}{1-\pi}m^{(k)}$, in which $p$ is the
point process defined by
\begin{equation}\label{eq:defn-p}
p((0,t]\times A) \triangleq \sum_{k=1}^{\infty}1_{\{\sigma_k \leq
t\}}1_{\{Y_k \in A\}}, \quad \quad t \geq 0, \, A \in
\mathcal{B}(\mathbb{R}^d).
\end{equation}

\end{lemma}

\begin{proof}
Using Bayes' formula, and the independence of the stochastic
elements $\theta$, $\Lambda$ and $X$ we can write
\begin{equation}\label{eq:phik-ik-vk}
\pk=\frac{\E_0\left[(\Lambda-\mu)^{k}Z_t 1_{\{\theta \leq t
\}}|\mathcal{F}_t\right]}{(1-\Pi_t) \E_0[Z_t|\mathcal{F}_t]} =
U_t^{(k)}+V_t^{(k)}
\end{equation}
in which
\begin{equation}\label{eq:defn-u}
U_t^{(k)} \triangleq \frac{\pi}{1-\pi}e^{\lambda t}
\int_{\mathbb{R}_+} (\nu-\mu)^k L^{\nu}_{t} \gamma(d \nu), \quad
\text{and}
\end{equation}
\begin{equation}\label{eq:defn-v}
V_t^{(k)}\triangleq \int_{0}^{t} \int_{\mathbb{R}_{+}}\lambda
e^{\lambda(t-u)} \frac{L^{\nu}_t}{L^{\nu}_u}(\nu-\mu)^{k}\gamma(d
\nu)du.
\end{equation}
Here we have used the notation
\begin{equation}\label{eq:defn-L-t}
L^{\nu}_t \triangleq e^{-(\nu-\mu)t}
\prod_{k=1}^{N_t}\left[\frac{\nu}{\mu} r(Y_k)\right], \quad \nu
\in \mathbb{R_+}.
\end{equation}
To derive (\ref{eq:phik-ik-vk}) we have used (\ref{eq:Z}),
(\ref{eq:defn-phik}) and the identity
\[
1-\Pi_t=\frac{(1-\pi)e^{-\lambda t}}{\E_0[Z_t|\mathcal{F}_t]},
\]
which we can derive using the independence of $\theta$ and $X$
under $\P_0$.

The process $L^{\nu}$ is the unique locally bounded solution of
the equation (see e.g. \cite{elliott})
\begin{equation}\label{eq:dynamics-of-L}
dL^{\nu}_t=L^{\nu}_{t -}\left[-(\nu-\mu)dt+\int_{y \in
\mathbb{R}^d}\left(\frac{\nu}{\mu}r(y)-1\right)p(dtdy)\right],
\end{equation}
with $L_0=1$. Using (\ref{eq:dynamics-of-L}) and the change of
variable formula it is easy to obtain
\begin{equation}
dU_t^{(k)}=(\lambda U_t^{(k)}-U_t^{(k+1)})dt+ \int_{y \in
\mathbb{R}^d}\left((r(y)-1)U_t^{(k)}+
\frac{r(y)}{\mu}U_t^{(k+1)}\right)p(dtdy),
\end{equation}
with $U^{(k)}_0=\frac{\pi}{1-\pi}m^{(k)}$, and
\begin{equation}
dV_t^{(k)}= (\lambda m^{(k)}-V_t^{(k+1)}+\lambda V_t^{(k)})dt +
\int_{y \in \mathbb{R}^d} \left((r(y)-1)V_t^{(k)}+
\frac{r(y)}{\mu}V_t^{(k+1)} \right)p(dtdy).
\end{equation}
with $V^{k}_0=0$. Now (\ref{eq:dyn-Phi}) follows from
(\ref{eq:phik-ik-vk}).
\end{proof}

Lemma~\ref{lem:dyn-of-Phi} shows that 1)$\Phi^{(0)}$ is not a
Markov process, and 2) the sequence $\{\Phi^{(k)}\}_{k \in
\mathbb{N}}$ has the Markovian property and its natural filtration
is the same as $\mathbb{F}$. The following corollary gives a
sufficient condition that the distribution of the post-disorder
rate $\gamma$ must satisfy in order for the sufficient statistic
to be finite dimensional.

\begin{cor}\label{cor:finite-dimensional}
If $\gamma$ is a discrete distribution with only $k$ atoms then
$\{\Phi^{(0)},\Phi^{(1)}, \cdots, \Phi^{(k-1)}\}$ is a
$k$-dimensional Markovian sufficient statistic.
\end{cor}
\begin{proof}
This follows from the same line of arguments used in the proof of
Corollary 3.3 in \cite{bdk05}. Here, we will give it not only for
readers conveneience but also because the notation we introduce
here will be used later. Let us denote by $\nu_1, \cdots, \nu_k$
the atoms of the distribution $\gamma$ and define
\begin{equation}
p(v) \triangleq \prod_{k=1}^{k}(v-\nu_i+\mu)\equiv
v^{k}+\sum_{i=0}^{k-1}c_i v^i, \quad v \in \mathbb{R},
\end{equation}
for some suitable numbers $c_0,...,c_{k-1}$. Observe that the
random variable
\[
p(\Lambda-\mu)=(\Lambda-\mu)^{k}+\sum_{i=0}^{k-1}c_i(\Lambda-\mu)^{i}=0,
\quad \text{a.s}.
\]
The last identity together with (\ref{eq:defn-phik}) implies that
\begin{equation}\label{eq:k-dimensional}
\Phi_{t}^{(k)}+\sum_{i=0}^{k-1}c_i \Phi_t^{(i)}=0, \quad \P-a.s.
\end{equation}
Now, it can be seen from the form of the penalty function in
(\ref{eq:new-penalty}) and the dynamics in (\ref{eq:dyn-Phi}) that
$\{\Phi^{(0)},\Phi^{(1)}, \cdots, \Phi^{(k-1)}\}$ is a
$k$-dimensional Markovian sufficient statistic.
\end{proof}

In the remainder of the paper we will assume that the distribution
for the post-disorder rate $\Lambda$ has Bernoulli distribution.

\section{Post-Disorder Rate with Bernoulli Distribution}\label{sec:Bernoulli}

In this section we will assume that the random variable $\Lambda$
takes either the value $\mu_1>0$ or $\mu_2>0$, i.e.
$\gamma(\{\mu_1,\mu_2\})=1$. From (\ref{eq:k-dimensional}) it
follows that $\Phi^{(2)}=(\mu_1+\mu_2-2
\mu)\Phi^{(1)}-(\mu_1-\mu)(\mu_2-\mu)\Phi^{(0)}$. According to
Lemma~\ref{lem:dyn-of-Phi}, the pair $(\Phi^{(0)},\Phi^{(1)})$
satisfies
\begin{equation}
\begin{split}
&
d\Phi^{(0)}_t=(\lambda(1+\Phi_t^{(0)})-\Phi_t^{(1)})dt+\Phi_{t-}^{(0)}\int_{y
\in \mathbb{R}^{d}}(r(y)-1)p(dt dy)+
\Phi_{t-}^{(1)}\frac{1}{\mu}\int_{y \in \mathbb{R}^{d}}r(y)p(dt
dy)
\\&d\Phi^{(1)}_t=(\lambda m^{(1)}+(\lambda-(\mu_1+\mu_2-2
\mu))\Phi^{(1)}_t+(\mu_1-\mu)(\mu_2-\mu)\Phi^{(0)}_t)dt
\\& +\Phi_{t-}^{(1)}\int_{y \in \mathbb{R}^{d}}(r(y)-1)p(dt dy)+
((\mu_1+\mu_2-2
\mu)\Phi^{(1)}_{t-}-(\mu_1-\mu)(\mu_2-\mu)\Phi_{t-}^{(0)})\frac{1}{\mu}\int_{y
\in \mathbb{R}^{d}}r(y)p(dt dy)
\end{split}
\end{equation}
with initial conditions $\Phi^{(0)}_0=\frac{\pi}{1-\pi}$ and
$\Phi^{(1)}_0=\frac{\pi}{1-\pi}m^{(1)}$.

 Instead of the sufficient statistic
$(\Phi^{(0)},\Phi^{(1)})$, it will be more convenient to work with
\begin{equation}\label{eq:sufficient-stats}
\tilde{\Phi}^{(0)}_t \triangleq \frac{\P\left\{\Lambda=\mu_1, \,
\theta \leq t|\mathcal{F}_t\right\}}{\P\{\theta>t|\mathcal{F}_t\}}
\quad \text{and} \quad  \tilde{\Phi}^{(1)}_t \triangleq
\frac{\P\left\{\Lambda=\mu_2,\, \theta \leq
t|\mathcal{F}_t\right\}}{\P\{\theta>t|\mathcal{F}_t\}}.
\end{equation}
In fact the following a one-to-one relationship between these two
pairs holds
\begin{equation}
\tilde{\Phi}_t^{(0)}=\frac{(\mu_2-\mu)\Phi_t^{(0)}-\Phi^{(1)}_t}{\mu_2-\mu_1}\,
 \quad \text{and} \quad
\tilde{\Phi}_t^{(1)}=\frac{(\mu_1-\mu)\Phi_t^{(0)}-\Phi^{(1)}_t}{\mu_1-\mu_2}\,
.
\end{equation}
The dynamics of this new sufficient statistic are autonomous as
can be seen from
\begin{equation}\label{eq:autonomous}
\begin{split}
d \tilde{\Phi}^{(0)}_{t}&=
\left\{\frac{\lambda(\mu_2-\mu-m^{(1)})}{\mu_2-\mu_1}+(\lambda-\mu_1+\mu)
\tilde{\Phi}_t^{(0)} \right\}dt+ \tilde{\Phi}_{t-}^{(0)}\int_{y
\in \mathbb{R}^d}\left[\left(1+\frac{\mu_1-\mu}{\mu
}\right)r(y)-1\right]p(dt dy)
\\ d \tilde{\Phi}^{(1)}_{t}&=
\left\{\frac{\lambda(\mu_1-\mu-m^{(1)})}{\mu_1-\mu_2}+(\lambda-\mu_2+\mu)
\tilde{\Phi}_t^{(1)} \right\}dt+  \tilde{\Phi}_{t-}^{(1)}\int_{y
\in \mathbb{R}^d}\left[\left(1+\frac{\mu_2-\mu}{\mu
}\right)r(y)-1\right]p(dt dy)
\end{split}
\end{equation}
When the number of atoms of the distribution $\gamma$ is more than
two, we expect that sufficient statistics defined similarly will
also be autonomous.

The sufficient statistic we introduced in
(\ref{eq:sufficient-stats}) has a natural interpretation and is
similar in flavor to particle filters: these are the normalized
probabilities that are assigned to each atom $\mu_i$ and these are
updated continuously between the times of the observations, since
not having an observation in fact reveals some information about
the intensity of the underlying Poisson process. Indeed from
(\ref{eq:autonomous}) we observe that the sufficient statistic
$(\tilde{\Phi}^{(0)},\tilde{\Phi}^{(1)})$ solves an ordinary
differential equation between the observations, and the terms that
involve the counting process $p$ are inactive. When there is an
observation, these normalized probabilities jump depending on the
jump size of the observation. We will see the optimal alarm
mechanism is to sound the alarm as soon as the sufficient
statistic touches or jumps above a convex and decreasing curve in
$\mathbb{R}_{+}^2$, if the sufficient statistic starts below this
curve. Otherwise it is optimal to sound the alarm immediately.
Since the jump distribution also changes at the time of disorder,
not only the timing of the observations but also the magnitude of
the observations is informative. Therefore, it is reasonable to
expect that we are able to construct a more acute alarm in this
case than the case in which the observations are coming from a
simple Poisson process where the jump size does not carry any
information.

In the case when the post disorder rate could go both up and down
by one unit, i.e., $\mu_1=\mu-1$ and $\mu_2=\mu+1$, then the
dynamics in (\ref{eq:autonomous}) become
\begin{equation}\label{eq:above-below}
\begin{split}
d \tilde{\Phi}^{(0)}_{t}&=
\left\{\frac{\lambda(1-m)}{2}+(\lambda+1) \tilde{\Phi}_t^{(0)}
\right\}dt+ \tilde{\Phi}_{t-}^{(0)}\int_{y \in
\mathbb{R}^d}\left[\left(1-\frac{1}{\mu }\right)r(y)-1\right]p(dt
dy)
\\ d \tilde{\Phi}^{(1)}_{t}&=
\left\{\frac{\lambda(1+m)}{2}+(\lambda-1) \tilde{\Phi}_t^{(1)}
\right\}dt+  \tilde{\Phi}_{t-}^{(1)}\int_{y \in
\mathbb{R}^d}\left[\left(1+\frac{1}{\mu }\right)r(y)-1\right]p(dt
dy),
\end{split}
\end{equation}
in which $m=m^{1}=\P\{\Lambda=\mu+1\}-\P\{\Lambda=\mu-1\} \in
[-1,1]$. Observe that when an arrival comes,
$\tilde{\Phi}^{(0)}_{t}$ jumps down and $\tilde{\Phi}^{(1)}$ jumps
up. Assuming $m \in (-1,1)$ then $\tilde{\Phi}^{(0)}_{t}$ is
always increasing between the observations. $\tilde{\Phi}^{(1)}$,
on the other hand, can be increasing or mean reverting depending
on the value of $\lambda$. Note that the values $m=-1$ or $m=1$
correspond to the degenerate cases in which the post-disorder rate
is known and the sufficient statistic becomes one-dimensional.

On the other hand, in the case when the post disorder rate could
only go up by one or two units, i.e., $\mu_1=\mu+1$ and
$\mu_2=\mu+2$, then the dynamics in (\ref{eq:autonomous}) become
\begin{equation}\label{eq:up-up}
\begin{split}
d \tilde{\Phi}^{(0)}_{t}&= \left\{\lambda(2-m)+(\lambda-1)
\tilde{\Phi}_t^{(0)} \right\}dt+ \tilde{\Phi}_{t-}^{(0)}\int_{y
\in \mathbb{R}^d}\left[\left(1+\frac{1}{\mu
}\right)r(y)-1\right]p(dt dy),
\\ d \tilde{\Phi}^{(1)}_{t}&=
\left\{\lambda(m-1)+(\lambda-2) \tilde{\Phi}_t^{(1)} \right\}dt+
\tilde{\Phi}_{t-}^{(1)}\int_{y \in
\mathbb{R}^d}\left[\left(1+\frac{2}{\mu }\right)r(y)-1\right]p(dt
dy),
\end{split}
\end{equation}
in which $m=2 \P\{\lambda=\mu+2\}+\P\{\lambda=\mu+1\} \in [1,2]$.
Here the initial conditions are
$\tilde{\Phi}^{(0)}_{0}=(2-m)\frac{\pi}{1-\pi}$ and
$\tilde{\Phi}^{(1)}_{0}=(m-1)\frac{\pi}{1-\pi}$. We will assume
that $m \in (1,2)$ as otherwise the problem degenerates into a
one-dimensional one. In the next section we will see that the
intuition that a decision would sound the alarm only at the times
when it observes an arrival does not in general hold; see
Remark~\ref{rem:only-at-jump-times}. This intuition becomes
relevant only when $\lambda$ and $c$ are small enough, i.e. when
the disorder intensity and delay penalty are small. If $\lambda
\geq 2$ then both $\tilde{\Phi}^{(0)}_{t}$ and
$\tilde{\Phi}^{(1)}_{t}$ increase between the jumps, because the
rate of disorder is high enough despite the fact that there have
been no arrivals. When $\lambda \in [1,2)$,
$\tilde{\Phi}^{(0)}_{t}$ increases between the jumps and
$\tilde{\Phi}^{(1)}_{t}$ is mean reverting. When $\lambda \in
(0,1)$, both $\tilde{\Phi}^{(0)}_{t}$ and $\tilde{\Phi}^{(1)}_{t}$
have mean reverting paths between arrivals. Since the post
disorder arrival rate can only move up, both
$\tilde{\Phi}^{(0)}_{t}$ and $\tilde{\Phi}^{(1)}_{t}$ have an
upward jump when there is an observation.

In the remainder of the paper we analyze the case when the
sufficient statistic is of the form (\ref{eq:up-up}). Note that in
this case the penalty function in (\ref{eq:new-penalty}) becomes
\begin{equation}\label{eq:penalty-up-up}
R_{\tau}(\pi)=1-\pi+
c(1-\pi)\mathbb{E}_0\left[\int_0^{\tau}e^{-\lambda
t}\left(\tilde{\Phi}_t^{(0)}+\tilde{\Phi}_t^{(1)}-\frac{\lambda}{c}\right)dt\right],
\quad \tau \in \mathcal{S}.
\end{equation}
Let us define
\begin{equation}\label{eq:explicit-x-and-y}
\begin{split}
 x(t,x_0)& \triangleq
\begin{cases}
-\frac{\lambda(2-m)}{\lambda-1}+e^{(\lambda-1)t}\left[x_0+\frac{\lambda(2-m)}{\lambda-1}\right],
& \lambda \neq  1, \\ x_0+(2-m)t, & \lambda=1,
\end{cases} \quad \text{and}
\\ y(t,y_0) & \triangleq
\begin{cases}
 -\frac{\lambda(m-1)}{\lambda-2}+
e^{(\lambda-2)t}\left[y_0+\frac{\lambda(m-1)}{\lambda-2}\right] &
\lambda \neq 2,
\\ y_0+2(m-1)t & \lambda=2.
\end{cases}
\end{split}
\end{equation}
Note that $x$ and $y$ satisfy the semigroup property, i.e., for
every $t \in \mathbb{R}$ and $s \in \mathbb{R}$,
\begin{equation}\label{eq:semi-group}
x(t+s,x_0)=x(s,x(t,x_0)) \quad \text{and} \quad
y(t+s,x_0)=y(s,y(t,x_0)).
\end{equation}
Let us denote by $\sigma_n$ the jump times of the process $X$.
Then we get
\begin{equation}\label{eq:paths-at-the-jumps}
\begin{split}
&
\tilde{\Phi}^{(0)}_{t}=x(t-\sigma_n,\tilde{\Phi}^{(0)}_{\sigma_n})
\quad \text{and} \quad
\tilde{\Phi}^{(1)}_{t}=y(t-\sigma_n,\tilde{\Phi}^{(1)}_{\sigma_n}),
\quad \sigma_n \leq t <\sigma_{n+1},
\\ &   \tilde{\Phi}^{(0)}_{\sigma_n}=
\left(1+\frac{1}{\mu}\right)r(Y_n) \tilde{\Phi}^{(0)}_{\sigma_n-}
\quad \text{and} \quad \tilde{\Phi}^{(1)}_{\sigma_n}=
\left(1+\frac{2}{\mu}\right)r(Y_n)\tilde{\Phi}^{(1)}_{\sigma_n-}
\quad  n \in \mathbb{N}_0.
\end{split}
\end{equation}
The minimum of the Bayes risk in (\ref{eq:penalty-up-up}) is given
by;
\begin{equation}
U(\pi)=\inf_{\tau \in
\mathcal{S}}R_{\tau}(\pi)=(1-\pi)+c(1-\pi)V\left((2-m)\frac{\pi}{1-\pi},(m-1)\frac{\pi}{1-\pi}\right),
\end{equation}
in which $V$ is defined as the value function of the optimal
stopping problem for a two-dimensional Markov process
\begin{equation}\label{eq:value-function}
V(\phi_0,\phi_1) \triangleq \inf_{\tau \in
\mathcal{S}}\E_0^{\phi_0,\phi_1}\left[\int_0^{\tau}e^{-\lambda t}
g\left( \tilde{\Phi}_t\right)dt\right], \quad \quad \tilde{\Phi}_t
\triangleq (\tilde{\Phi}^{(0)}_{t},\tilde{\Phi}^{(1)}_{t}),
\end{equation}
with a running cost function
\begin{equation}
g(\phi_0,\phi_1)=\phi_0+\phi_1-\frac{\lambda}{c}.
\end{equation}
Here, $\mathbb{E}_0^{\phi_0,\phi_1}$ is the conditional $\P_0$
expectation given that $\tilde{\Phi}^{(0)}_{0}=\phi_0$ and
$\tilde{\Phi}^{(1)}_{0}=\phi_1$.

\section{Upper and Lower Bounds on the Optimal Stopping
Time}\label{sec:bounds}

Unlike the optimal stopping problem for It\^{o} diffusions,
analyzing the sample path behavior of the piece-wise deterministic
Markov process $\tilde{\Phi} \triangleq (\tilde{\Phi}_1,
\tilde{\Phi}_2)$, we are able to determine the optimal stopping
time for most parameter values. For remaining parameter values we
are able to provide some lower bound and an upper bounds on the
optimal stopping time.

All the results in this section assume that an optimal stopping
time exists and it is given by
\begin{equation}\label{optst}
\tau^{*}(\phi_0,\phi_1) \triangleq \inf\{t \geq 0: V(\p_t)=0,\,
\p_0=(\phi_0,\phi_1)\}.
\end{equation}
In Section~\ref{sec:optimal-stopping-time}, we verify that this
assumption in fact holds. With (\ref{optst}) we will call the
region
\begin{equation}
\bGamma \triangleq \{(\phi_0,\phi_1)\in \R^2_+:
  v(\phi_0,\phi_1)=0\},  \quad \C \triangleq \R^2_+ \backslash
  \bGamma,
\end{equation}
the \emph{optimal stopping region}. Let us start this section with
a simple observation.
\begin{lemma}
Let us define
\begin{equation}\label{defn:tau-l}
\tau^{l} \triangleq \inf\{t\geq 0: \tilde{\Phi}^{(0)}_t +
\tilde{\Phi}^{(1)}_t \geq \lambda/c\}.
\end{equation}
If there is an optimal stopping time for the problem in
(\ref{eq:value-function}), let us denote it by $\tau^*$, then
$\tau^* \geq \tau^l$.
\end{lemma}
\begin{proof}
Let $\tau \in \mathcal{S}$ be any stopping rule. Then
\begin{equation}
\begin{split}
\E_0^{\phi_0,\phi_1}\left[\int_0^{\tau \vee \tau^l}e^{-\lambda t
}g(\tilde{\Phi}_t)dt\right] &=
\E_0^{\phi_0,\phi_1}\left[\int_0^{\tau}e^{-\lambda t
}g(\tilde{\Phi}_t)dt\right]+\E_0^{\phi_0,\phi_1}\left[1_{\{\tau^l>\tau\}}\int_{\tau}^{\tau^l}e^{-\lambda
t }g(\tilde{\Phi}_t)dt\right]
\\ & \leq \E_0^{\phi_0,\phi_1}\left[\int_0^{\tau}e^{-\lambda t
}g(\tilde{\Phi}_t)dt\right], \quad (\phi_0,\phi_1) \in
\mathbb{R}_+^2.
\end{split}
\end{equation}
Here $\tau \vee \tau^l=\max\{\tau, \tau^l\}$.
\end{proof}

When the rate of disorder or $c$ in (\ref{eq:penalty}) are large
enough, then in fact the lower bound $\tau^l$ is optimal as the
following proposition illustrates, i.e. the free boundary
corresponding to the two-dimensional optimal stopping problem in
(\ref{eq:value-function}) can be determined completely. This is a
very special instance of a multi-dimensional optimal stopping
problem where an explicit determination of the free boundary is
possible.
\begin{proposition}\label{prop:explicit-solution}
If (i) $\lambda \geq 2$, or, (ii) $\lambda \in [1,2)$ and $c \geq
2-\lambda$, or, (iii) $\lambda \in (0,1)$ and $c \geq \max
\left(2-\lambda,1-\lambda\right)$, then the stopping rule $\tau^l$
of (\ref{defn:tau-l}) is optimal for the problem in
(\ref{eq:value-function}).
\end{proposition}

\begin{proof}
(i) Let us first consider the case $\lambda \geq 2$. It is clear
from the dynamics of the sufficient statistic in (\ref{eq:up-up})
that the sample paths of $\tilde{\Phi}^{(0)}_t$ and
$\tilde{\Phi}^{(1)}_t$ are increasing functions of time. Therefore
the process $\tilde{\Phi}$ does not return to the region
$\{(\phi_0,\phi_1) \in \mathbb{R}^2_+: \phi_0+\phi_1 \leq
\lambda/c\}$. Thus for every stopping time $\tau \in \mathcal{S}$
\begin{equation}
\begin{split}
&\E_0^{\phi_0,\phi_1} \left[\int_0^{\tau}e^{-\lambda t
}g(\tilde{\Phi}_t)dt\right] \geq
\E_0^{\phi_0,\phi_1}\left[\int_0^{\tau \vee \tau^l}e^{-\lambda t
}g(\tilde{\Phi}_t)dt\right] \\ &=
\E_0^{\phi_0,\phi_1}\left[\int_0^{ \tau^l}e^{-\lambda t
}g(\tilde{\Phi}_t)dt\right]+ \E_0^{\phi_0,\phi_1}\left[1_{\{\tau
\geq \tau^l\}}\int_{\tau^l}^{ \tau}e^{-\lambda t
}g(\tilde{\Phi}_t)dt\right] \geq
\E_0^{\phi_0,\phi_1}\left[\int_0^{ \tau^l}e^{-\lambda t
}g(\tilde{\Phi}_t)dt\right]
\end{split}
\end{equation}
(ii) If $\lambda \in [1,2)$ then any sample path of
$\tilde{\Phi}^{(0)}$ is still an increasing function of $t$, but
the same is not true anymore for the sample paths of
$\tilde{\Phi}^{(1)}$. The paths of $\py$ increase with jumps;
between the jumps the paths are mean reverting to the level
$\lambda(m-1)/(2-\lambda)$. However, since the processes $\px$ and
$\py$ can only increase by jumps we have that
\begin{equation}
\px_t \geq x(t,\phi_0) \quad \text{and} \quad \py_t \geq
y(t,\phi_1), \quad t \geq 0.
\end{equation}
Therefore
\begin{equation}\label{eq:V-vs-deterministic}
 V(\phi_0,\phi_1) \geq \inf_{\tau \in \mathcal{S}
}\E_0^{\phi_0,\phi_1}\left[\int_{0}^{\tau}e^{-\lambda t}
\left(x(t,\phi_0)+y(t,\phi_1)-\frac{\lambda}{c}\right) dt\right].
\end{equation}
Clearly if for any $(\phi_0,\phi_1)$ if the right hand side of
(\ref{eq:V-vs-deterministic}) is zero, then $V=0$, since we also
know that $V \leq 0$. This can be used to find a superset of the
continuation region. However, as we shall see shortly this
superset coincides with the \emph{advantageous region}
\begin{equation}\label{defn:advantageous}
\mathbb{C}_0 \triangleq \{(\phi_0,\phi_1) \in \mathbb{R}_+^2:
\phi_0+\phi_1 \leq \lambda/c\}.
\end{equation}
Observe that it is not optimal to stop before $\tilde{\Phi}$
leaves the region $\mathbb{C}_0$.

Let us take a look at the derivative of the integrand on the
righthand side in (\ref{eq:V-vs-deterministic}),
\begin{equation}\label{eq:derivative-sum}
\frac{d}{dt}[x(t,\phi_0)+y(t,\phi_1)]=(\lambda-1)x(t,\phi_0)+(\lambda-2)y(t,\phi_1)+\lambda.
\end{equation}
The righthand side of (\ref{eq:derivative-sum}) vanishes if the
curve $t \rightarrow (x(t,\phi_0), y(t,\phi_1))$ meets the line
\begin{equation}\label{eq:line}
l: (\lambda-1)x+(\lambda-2)y+\lambda=0.
\end{equation}
Note that since $\lambda \in [1,2)$ the y-intercept of the line is
such that $\frac{\lambda}{2-\lambda} \geq \lambda$. Since $l$ is
increasing and $c \geq 2-\lambda$, the intersection of $l$ with
the set $\mathbb{C}_0$ is empty. Observe also that every $t
\rightarrow (x(t,\phi_0),y(t,\phi_1))$ starting at
$(\phi_0,\phi_1)$ is decreasing and the derivative in
(\ref{eq:derivative-sum}) is increasing. Therefore, $t \rightarrow
(x(t,\phi_0),y(t,\phi_1))$ meets the line $l$ at most once for any
$(\phi_0,\phi_1) \in \mathbb{R}_{+}$.

\begin{equation}\label{eq:text}
\begin{cases}
\text{ Furthermore, if $t \rightarrow (x(t,\phi_0),y(t,\phi_1))$
meets $l$ at $t_{l}=t_{l}(\phi_0,\phi_1)$, then the function} \\
\text{$t \rightarrow x(t,\phi_0)+y(t,\phi_1)$ is decreasing on
$[0,t_l]$ and increasing on $[t_l,\infty)$. If $t \rightarrow
 (x(t,\phi_0)$} \\ \text{ $y(t,\phi_1))$  does not intersect $l$, then the
function $t \rightarrow x(t,\phi_0)+y(t,\phi_1)$ is increasing on}
\\ [0,\infty).
\end{cases}
\end{equation}

Since the line $l$ does not meet the region $\mathbb{C}_0$ for
every $(\phi_0,\phi_1) \in l$ we have that $\phi_0+\phi_1 \geq
\lambda/c$. Now (\ref{eq:text}) implies that
$x(t,\phi_0)+y(t,\phi_1)-\frac{\lambda}{c}>0$ for $(\phi_0,\phi_1)
\in \mathbb{R}_{+} ^2-\mathbb{C}_0$ and $t \geq 0$. This implies
that the righthand side of (\ref{eq:V-vs-deterministic}) is zero,
which in turn implies that $V(\phi_0,\phi_1)=0$ for all
$(\phi_0,\phi_1) \in \mathbb{R}_{+} ^2-\mathbb{C}_0$.

(iii) If $\lambda \in (0,1)$, then both of the paths of
$x(t,\phi_0)$ and $y(t,\phi_1)$ are mean reverting. Because of our
assumption on $c$ the line $l$ in (\ref{eq:line}) does not
intersect with $\mathbb{C}_0$ and lies entirely above this region.
Let us denote the region between $l$ and $\mathbb{C}_0$ by
\begin{equation}
\text{Sh}\triangleq\{(\phi_0,\phi_1) \in \R_{+}^2:
\phi_0+\phi_1-\lambda/c>0, (\lambda-1)\phi_0+ (\lambda-2)
\phi_1+\lambda< 0\}.
\end{equation}
From (\ref{eq:derivative-sum}) it follows that
$x(t,\phi_0)+y(t,\phi_1)>\lambda/c$ if $(\phi_0,\phi_1) \in
\text{Sh}$. Therefore, the path $t\rightarrow
(x(t,\phi_0),y(t,\phi_1))$ never enters the region $\mathbb{C}_0$
if $(\phi_0,\phi_1) \notin \mathbb{C}_0$. Therefore, the righthand
side of (\ref{eq:V-vs-deterministic}) is zero, which in turn
implies that $V(\phi_0,\phi_1)=0$ for any $(\phi_0,\phi_1) \in
\mathbb{R}^{2}_+-\mathbb{C}_0$.
\end{proof}

\begin{proposition}\label{prop:lambda-1-2}
Assume $\lambda \in [1,2)$ and $c \in (0,2-\lambda)$. Let us
define
\begin{equation}
D \triangleq \{(\phi_0,\phi_1) \in \mathbb{R}_{+}^2: \phi_0 \leq
\phi_0^*,\,\,  \phi_0+\phi_1 \leq \xi\} \cup \{(\phi_0,\phi_1) \in
\mathbb{R}_{+}^2: \phi_0 >\phi_0^*,\,\,\phi_0+\phi_1 \leq
\lambda/c\},
\end{equation}
in which $(\phi^*_0,\phi^*_1) \triangleq
(\lambda(-1+(2-\lambda)/c),\lambda(1+(\lambda-1)/c))$ and
\[
\xi=y\left(-t^*,\lambda\left(\frac{\lambda-1}{c}+1\right)\right),
\quad \text{where} \quad x\left(-t^*,\lambda\left(\frac{
2-\lambda}{c}-1\right)\right)=0.
\]
Then the region $D$ is a superset of the optimal stopping region.
\end{proposition}

\begin{proof}
Let us note that (\ref{eq:V-vs-deterministic}) implies that
\begin{equation}\label{eq:V-geq-deterministic}
 V(\phi_0,\phi_1) \geq \inf_{t \in [0,\infty]
}\left[\int_{0}^{t}e^{-\lambda s}
\left(x(s,\phi_0)+y(s,\phi_1)-\frac{\lambda}{c}\right) ds\right].
\end{equation}
Because of the assumption on $c$ the line in (\ref{eq:line})
intersects the region $\mathbb{C}_0$ defined in
(\ref{defn:advantageous}). Note that $l$ and the boundary
$x+y-\lambda/c=0$ of the region $\mathbb{C}_0$ intersect at
$(\phi^*_0,\phi^*_1)$. By running the time ``backwards'', we can
find $\xi$ and $t^*$ such that
\begin{equation}
(0,\xi)=(x(-t^*,\phi_0^*),y(-t^*,\phi_1^*)).
\end{equation}
By the semi-group property (see (\ref{eq:semi-group})), we have
\[
x(t^*,0)=x(t^*,x(-t^*,\phi_0^*))=x(t^*+(-t^*),\phi_0^*)=x(0,\phi_0^*)=\phi_0^*,
\]
and,
\[
y(t^*,\xi)=y(t^*,x(-t^*,\phi_1^*))=y(t^*+(-t^*),\phi_1^*)=y(0,\phi_1^*)=\phi_1^*.
\]
So, the curve $t \rightarrow (x(t,0),y(t,\xi))$, $t \geq 0$, meets
line $l$ at $(\phi_0^*,\phi_1^*)$, and $t_l$ in (\ref{eq:text})
equals to $t^*$. This implies that
\[
x(t,0)+y(t,\xi) \geq x(t^*,0)+y(t^*,\xi)=
\phi_0^*+\phi_1^*=\frac{\lambda}{c},
\]
and in particular $\xi \geq \lambda/c$. Now we will show that when
$\lambda$ and $c$ are chosen as in the statement of the
proposition it is optimal to stop outside the region $D$.

The curve $t \rightarrow (x(t,0),y(t,\xi))$ divides
$\mathbb{R}_{+}^2$ into two connected components containing
$\mathbb{C}_0$ and the region
\begin{equation}
M \triangleq \mathbb{R}_{+}^2-D) \cap \{(x,y) \in
\mathbb{R}_{+}^2: (\lambda-1)x+(\lambda-2)y +\lambda<0 \}
\end{equation}
respectively. Every curve $t \rightarrow
(x(t,\phi_0),y(t,\phi_1))$, $t \geq 0$ starting at
$(\phi_0,\phi_1) \in M$ will stay in $M$, since from the
semi-group property (\ref{eq:semi-group}) it follows that two
distinct curves $t\rightarrow (x(t,\phi^a_0),y(t,\phi_1^a))$ and
$t\rightarrow (x(t,\phi^b_0),y(t,\phi_1^b))$ do not intersect.
Therefore, $t \rightarrow (x(t,\phi_0),y(t,\phi_1))$, $t \geq 0$,
$(\phi_0,\phi_1) \in M$ intersects the line $l$ in (\ref{eq:line})
away from $\mathbb{C}_0$ and (\ref{eq:text}) implies that
$x(t,\phi_0)+y(t,\phi_1) > \lambda/c$ for any $(\phi_0,\phi_1) \in
M$. Now, from (\ref{eq:V-geq-deterministic}) we conclude that
$V=0$ since the infimum on the right-hand-side is equal to 0 from
the arguments above and we already know that $V \leq 0$.

On the other hand, if $(\phi_0,\phi_1) \in (\mathbb{R}_{+}^2-D)
\cap \{(x,y) \in \mathbb{R}_{+}^2: (\lambda-1)x+(\lambda-2)y
+\lambda \geq 0 \}$, the curve $t \rightarrow
(x(t,\phi_0),y(t,\phi_1))$, $t \geq 0$ does not intersect the line
$l$; therefore, the function $t \rightarrow
x(t,\phi_0)+y(t,\phi_1)$ is increasing and
\[
x(t,\phi_0)+y(t,\phi_1)>x(0,\phi_0)+y(0,\phi_1)\geq \phi_0+\phi_1
\geq \xi \geq \frac{\lambda}{c}, \quad 0<t<\infty.
\]
Again, the infimum on the right-hand-side of
(\ref{eq:V-geq-deterministic}) is equal to zero, which implies
that $V=0$.

\end{proof}

\begin{remark}\label{rem:only-at-jump-times}
If $\lambda \in [1,2)$ and $c \in (0,2-\lambda)$, then the
following line segment is a subset of the free boundary
\begin{equation}
H \triangleq \left\{(\phi_0,\phi_1) \in
\mathbb{R}_+^2:\phi_0+\phi_1 -\frac{\lambda}{c} = 0,\,\,\, \phi_1
\leq \phi^*_1 \right\}.
\end{equation}
This set in fact in the entrance boundary of the stopping region
(the boundary through which the path $t \rightarrow
(x(t,\phi_0),y(t,\phi_1))$ enters the stopping region).

\end{remark}

\begin{proposition}\label{prop:stop-at-jumps}
Assume that $\lambda \in (0,1)$ and that $0< c \leq
\frac{(2-\lambda)(1-\lambda)}{3-\lambda-m}$. If furthermore $c
\geq 2 \frac{1-\lambda}{3-m}$, then
\begin{equation}
P \triangleq \left\{(\phi_0,\phi_1) \in \mathbb{R}_+^2:
\phi_0+\frac{1}{2}\phi_1+\frac{3}{2}-\frac{1}{2}m-\frac{1}{c}\geq
0, \,\,\, \phi_0+\phi_1 -\frac{\lambda}{c} \geq 0 \right\},
\end{equation}
 is a subset of the optimal stopping region.

If on the other hand, $0<c<2 \frac{1-\lambda}{3-m}$, then the
first time
 time $(\px,\py)$
reaches the set,
\begin{equation}
R \triangleq \left\{(\phi_0,\phi_1) \in \mathbb{R}_+^2:
\phi_0+\frac{1}{2}\phi_1+\frac{3}{2}-\frac{1}{2}m-\frac{1}{c}\geq
0\right\},
\end{equation}
 is an upper bound on the
optimal stopping time.
\end{proposition}

\begin{proof}
When $0< c \leq \frac{(2-\lambda)(1-\lambda)}{3-\lambda-m}$, then
the line $l \cap \R_+^2$ lies entirely in $\Co$. The paths, $t
\rightarrow x(t,\phi_0,\phi_1)$, $t \geq 0$, that do not originate
in $\mathbb{C}_0$ enter into this region through the boundary
$\{(\phi_0,\phi_1) \in \R_{+}^2: \phi_0+\phi_1=\lambda/c\}$ and
once they cross into $\mathbb{C}_0$ they never leave it again
since $x(t,\phi_0)+y(t,\phi_0)<\phi_0+\phi_1<\lambda/c$ for any
point $(\phi_0,\phi_1) \in \mathbb{C}_0 \cap \{(\phi_0,\phi_1) \in
\R_{+}^2: (\lambda-1)\phi_0+ (\lambda-2)\phi_1+\lambda<0\}$, which
follows from (\ref{eq:derivative-sum}). Therefore the infimum on
the right-hand-side of (\ref{eq:V-geq-deterministic}) is attained
by either $t=0$ or $t=\infty$ if $(\phi_0,\phi_1) \in
\mathbb{R}^{2}_+-\mathbb{C}_0$. Either one never stops, pays a
penalty for being outside $\mathbb{C}_0$ for a while and then
enjoys being in this region ad infinitum, or stops immediately
because the cost of the initial penalty is deterrent enough. Since
\begin{equation}
\begin{split}
&\int_{0}^{\infty}\left(x(t,\phi_0)+y(t,\phi_1)-\frac{\lambda}{c}\right)dt=\phi_0+\frac{1}{2}\phi_1+\frac{3}{2}-\frac{1}{2}m-
\frac{1}{c},
\end{split}
\end{equation}
the infimum on the right-hand-side of
(\ref{eq:V-geq-deterministic}) is attained by $t=0$ if
$(\phi_0,\phi_1) \in P$, which in turn implies that
$V(\phi_0,\phi_1)=0$ for any $(\phi_0,\phi_1) \in P$. Observe
that, if $0<c<2 \frac{1-\lambda}{3-m}$ then $P=R$.
\end{proof}

\begin{remark}
Observe that if $\lambda \in (0,1)$ and $2 \frac{1-\lambda}{3-m}
\leq c \leq \frac{(2-\lambda)(1-\lambda)}{3-\lambda-m}$, then the
following line segment is a subset of the free boundary;
\begin{equation}
F \triangleq \left\{(\phi_0,\phi_1) \in
\mathbb{R}_+^2:\phi_0+\phi_1 -\frac{\lambda}{c} = 0,\,\,\, \phi_1
\leq 2 \left(-\frac{1-\lambda}{c}+\frac{3-m}{2}\right) \right\}.
\end{equation}
This region is in the \emph{exit boundary} of the stopping region
(i.e., the boundary through which the path $t \rightarrow
(x(t,\phi_0),y(t,\phi_1))$ exits from the stopping region).
\end{remark}

\begin{proposition}\label{prop:sec5-last}
Assume that $\lambda \in (0,1)$ and that
\begin{equation}\label{eq:c-range}
\frac{(2-\lambda)(1-\lambda)}{3-\lambda-m} < c < \max
\left(2-\lambda,1-\lambda \right).
\end{equation}
Then the region $D$ defined in Proposition~\ref{prop:lambda-1-2}
is a superset of the optimal stopping region.

\end{proposition}

\begin{proof}
From the assumption on the parameters $\lambda$ and $c$ it follows
that the mean reversion level $M=
\left(\frac{\lambda(2-m)}{1-\lambda},\frac{\lambda(m-1)}{2-\lambda}\right)$
of the path $t \rightarrow (x(t,\phi_0),y(t,\phi_1))$, $t \geq 0$,
is in  the region $[0,\lambda/c]\times[0,\lambda/c]
-\mathbb{C}_0$. Also, one can easily check that $M \in l$, in
which $l$ is as in (\ref{eq:line}). Line $l$ and the boundary of
the region $\mathbb{C}_0$ intersect at $(\phi_0^*,\phi_1^*)$.
Because $c>(2-\lambda)(1-\lambda)/(3-\lambda-m)$, the equation (as
an equation in the $t$-variable) $x(t,0)=\phi_0^*$ has a positive
solution, $t^*$ and $y_0=y(-t^*,\phi_1^*)>0$. The rest of the
proof follows by using the same arguments as in the proof of
Proposition~\ref{prop:lambda-1-2}.

%
%
\end{proof}

 \begin{remark}
Observe that if $\lambda \in (0,1)$ and $c$ satisfies
(\ref{eq:c-range}), then the following line segment is a subset of
the free boundary
\begin{equation}
A \triangleq \left\{(\phi_0,\phi_1) \in
\mathbb{R}_+^2:\phi_0+\phi_1 -\frac{\lambda}{c} = 0,\,\,\, \phi_1
\leq \lambda \left(1-\frac{1-\lambda}{c}\right) \right\}.
\end{equation}
Moreover, this set is a subset of entrance boundary of the
stopping region.
\end{remark}

\begin{remark}
If the assumptions of Proposition~\ref{prop:stop-at-jumps} are
satisfied, then it is optimal to sound the alarm only at arrival
times of the observation. This corresponds to the case when the
mean reversion level of the paths $t \rightarrow
(x(t,\phi_0),y(t,\phi_1))$ is inside the advantageous region
$\mathbb{C}_0$, which is defined in (\ref{defn:advantageous}).
Otherwise, since the paths of the sufficient statistic, $t
\rightarrow \p_t$, may reach the stopping region continuously or
via jumps, it might be optimal to declare the alarm between two
observations.

\end{remark}

We will close this section by proving that the optimal stopping
time $\tau^*$ is finite almost surely.
\begin{proposition}
 Let $\eta$ be a positive number such that
the region $\{(\phi_0,\phi_1): \phi_0+\phi_1 \geq \eta\}$  is a
subset of the stopping region. (The existence of $\eta$ is
guaranteed by
Propositions~~\ref{prop:explicit-solution}-\ref{prop:sec5-last}).
Let us denote the hitting time of this region by $\tau^u$. Then
$\Ei_0[\tau^u] \leq\eta (2+1/\mu)$. This implies that $\tau^*$ is
finite $\P_0$ almost surely.
\end{proposition}

\begin{proof}
Since the compensator of $p(dt dy)$ (defined in (\ref{eq:defn-p}))
is equal to $ \mu \beta_0(y) $ we can write the dynamics of $\px$
in (\ref{eq:up-up}) as
\begin{equation}
\begin{split}
\px_{t \wedge \tau^u}&=\px_0+\int_{0}^{t \wedge \tau^u}
 \left\{\lambda(2-m)+(\lambda-1)
\tilde{\Phi}_t^{(0)} \right\}dt+ \int_0^{t \wedge \tau^u}\mu
\tilde{\Phi}_{t-}^{(0)}\int_{y \in
\mathbb{R}^d}\left[\left(1+\frac{1}{\mu
}\right)r(y)-1\right]\beta_0 (dy)ds
\\& +\int_0^{t \wedge \tau^u}\tilde{\Phi}_{t-}^{(0)}\int_{y
\in \mathbb{R}^d}\left[\left(1+\frac{1}{\mu
}\right)r(y)-1\right]q(ds dy)
\\&=\px_0+\int_{0}^{t \wedge \tau^u}
 \left\{\lambda(2-m)+\lambda
\tilde{\Phi}_t^{(0)} \right\}ds+\int_0^{t \wedge
\tau^u}\tilde{\Phi}_{t-}^{(0)}\int_{y \in
\mathbb{R}^d}\left[\left(1+\frac{1}{\mu }\right)r(y)-1\right]q(ds
dy),
\end{split}
\end{equation}
in which $q(dt dy) \triangleq p(dt dy)- \mu \beta_0(y) $ Here, we
have used the fact that $\int_{y \in \R_{+}^d}r(y)\beta_0(y)=1$.
The integral with respect to $q(dt dy)$ is an $\mathbb{F}$
martingale under the measure $\mathbb{P}_0$, since
\[
\begin{split}
\Ei\left[  \int_0^{t \wedge \tau^u}\mu
\tilde{\Phi}_{t-}^{(0)}\int_{y \in
\mathbb{R}^d}\left|\left(1+\frac{1}{\mu
}\right)r(y)-1\right|\beta_0 (dy)ds  \right] &\leq
\Ei_0\left[\int_0^{t \wedge
\tau^u}\left(2+\frac{1}{\mu}\right)\px_{s-}ds\right]
\\& \leq t \left(2+\frac{1}{\mu} \eta\right).
\end{split}
\]
Therefore
\[
\Ei_0\left[\px_{t \wedge
\tau^u}\right]=\phi_0+\Ei_0\left[\int_{0}^{t \wedge \tau^u}
 \left\{\lambda(2-m)+\lambda
\tilde{\Phi}_t^{(0)} \right\}ds\right] \geq
\lambda(2-m)\Ei_0\left[t \wedge \tau^u\right].
\]
On the other hand,
\[
\px_{t \wedge \tau^u} \leq \max\left(\eta,
\left(1+\frac{1}{\mu}\right)r(Y_{N_{t \wedge \tau^u}})\px_{t
\wedge \tau^u-}\right)\leq
\eta\left(1+\left(1+\frac{1}{\mu}\right)r(Y_{N_{t \wedge
\tau^u}})\right),
\]
almost surely; therefore
\[
\begin{split}
\Ei_0\left[t \wedge \tau^u\right] \leq \frac{1}{\lambda(2-m)}
\Ei_0\left[\px_{t \wedge \tau^u}\right] &\leq
\Ei_0\left[\eta\left(1+\left(1+\frac{1}{\mu}\right)r(Y_{1})\right)\right]
\\&=\eta \left(2+\frac{1}{\mu}\right).
\end{split}
\]
The result follows after an application of the monotone
convergence theorem.
\end{proof}

In what follows we will consider the cases in which the parameters
do not satisfy  the hypothesis of
Proposition~\ref{prop:explicit-solution} and construct a sequence
of functions iteratively, using an appropriately defined
functional operator, that converges to the value function
exponentially fast.

\section{Optimal Stopping Time and Properties of the Value Function and the Stopping Boundary}\label{sec:optimal-stopping-time} 
 The usual starting point to calculate the value function
in ( \ref{eq:value-function}) and find the optimal stopping time
would be to try to characterize the value function as the unique
solution of the free boundary problem
\begin{equation}\label{eq:free-boundary}
\min\{(\mathcal{A}-\lambda)v(\varphi)+g(\varphi),-v(\varphi)\}=0,
\end{equation}
in which the differential operator is the inifinitesimal generator
of the Markov process $(\tilde{\Phi}^{(0)},\tilde{\Phi}^{(1)})$
and whose action on a test function $f$ is given by
\begin{equation}
\begin{split}
\mathcal{A}f(\phi_0,\phi_1)&=\frac{\partial f}{\partial
\phi_0}(\phi_0,\phi_1)\left[\lambda(2-m)+(\lambda-1)
\phi_0\right]+\frac{\partial f}{\partial
\phi_1}(\phi_0,\phi_1)\left[\lambda(m-1)+(\lambda-2) \phi_1\right]
\\&+ \mu \int_{y \in
\mathbb{R}^d}\left[f\left(\left(1+\frac{1}{\mu}\right)r(y)\phi_0,
\left(1+\frac{2}{\mu}\right)r(y)\phi_1
\right)-f(\phi_0,\phi_1)\right] \beta_0(dy).
\end{split}
\end{equation}
The solution of the free boundary problem (\ref{eq:free-boundary})
may be identified by using certain boundary conditions (the smooth
fit principle). The smooth fit is expected to fail for
(\ref{eq:free-boundary}) at the exit boundary of the stopping
region. See e.g. \cite{BD04}, \cite{BD03} for failure of the
smooth fit principle when the infinitesimal generator
$\mathcal{A}$ is a differential delay operator (these papers
consider one dimensional free boundary problems). Instead of the
characterization of the value function as a solution of
quasi-variational inequalities, we will use a new characterization
of the value function of the optimal stopping problem in
(\ref{eq:value-function}). Specifically, we will  construct a
sequence of functions iteratively, using an appropriately defined
functional operator, that converges to the value function
exponentially fast. This will let us show that $\tau^*$ in
(\ref{optst}) is the optimal stopping time. We will also be able
to show the concavity of the value function and the convexity of
the free boundary.

\subsection{Optimal Stopping with Time Horizon
$\sigma_n$}

In this section, we will approximate the value function $V$ with a
sequence of optimal stopping problems. Let us denote
\begin{equation}\label{defnofVn}
V_{n}(\phi_0,\phi_1) \triangleq \inf_{\tau \in
\mathcal{S}}\E^{\phi_0,\phi_1}_0 \left[\int_0^{\tau \wedge
\sigma_n }e^{-\lambda t}g\left(\px_t,\py_t\right)dt \right]
\end{equation}
where $(\phi_0,\phi_1) \in \mathbb{R}^2_+$, $n \in \mathbb{N}$,
and $\sigma_n$ is the $n^{\text{th}}$ jump time of the process
$X$. Observe that $(V_n)_{n \in \mathbb{N}}$ is decreasing and
satisfies $-1/c<V_{n}<0$. Since $(\sigma_n)_{n \geq 1}$ is an
almost surely increasing sequence, $(V_{n})_{n \geq 1}$ is
decreasing. Therefore $\lim_{n}V_{n}$ exists. It is also immediate
that $V_n \geq V$. In fact we can say more about the limit of the
sequence $(V_{n})_{n \geq 1}$ as the next proposition illustrates.
\begin{proposition}\label{expofast}
 $V_{n}(\phi_0,\phi_1)$ converges to $V$ uniformly in
 $(\phi_0,\phi_1)\in \mathbb{R}_{+}^{2}$. In fact the rate of
 convergence is exponential as the following equation illustrates:
\begin{equation}
-\frac{1}{c}\left(\frac{\mu}{\mu+\lambda}\right)^n \geq
V_{n}(\phi_0,\phi_1)-V(\phi_0,\phi_1) \geq 0.
\end{equation}
\end{proposition}
\begin{proof}
\begin{equation}\label{eq:bndd}
\begin{split}
 \E^{\phi_0,\phi_1}_0 \left[\int_0^{\tau}e^{-\lambda
t}g\left(\tilde{\Phi}_t\right)dt \right] =\E^{\phi_0,\phi_1}_0
\left[\int_0^{\tau \wedge \sigma_n }e^{-\lambda
t}g\left(\tilde{\Phi}_t\right)dt \right] +
\E^{\phi_0,\phi_1}_0\left[1_{\{\tau \geq \sigma_n\}}
\int_{\sigma_n}^{\tau} e^{-\lambda t} g(\tilde{\Phi}_t)dt\right]
\end{split}
\end{equation}
The first term on the right-hand-side of (\ref{eq:bndd}) is
greater than $V_n$. Since $g(\cdot,\cdot)>-\lambda/c$ we can show
that the second term is greater than
\begin{equation}
-\frac{\lambda}{c}\E^{\phi_0,\phi_1}_0\left[1_{\{\tau \geq
\sigma_n\}}\int_{\sigma_n}^{\tau}e^{-\lambda s}ds\right] \geq
-\frac{1}{c} \E^{\phi_0,\phi_1}_0\left[e^{-\lambda
\sigma_n}\right] \geq -\frac{1}{c}\left(\frac{\mu}{\lambda+ \mu
}\right)^n.
\end{equation}
To show the last inequality we have used the fact that $\sigma_n$
is a sum of $n$ independent and identically distributed
exponential random variables with rate $\mu$ (i.e. $\sigma_n$ has
the Erlang distribution).
\end{proof}

Next, we will show that $V_n$ can be determined using an iterative
algorithm. To this end we introduce the following operators acting
on bounded Borel functions $f: \mathbb{R}_{+}^2 \rightarrow
\mathbb{R}$
\begin{align}
J f(t,\phi_0,\phi_1) & \triangleq
\E^{\phi_0,\phi_1}_{0}\bigg[\int_0^{t \wedge \sigma_1}e^{-\lambda
s} g(\px_s,\py_s)ds +1_{\{t \geq \sigma_1\}}e^{-\lambda \sigma_1}
f(\px_{\sigma_1},\py_{\sigma_1}) \bigg],\,\, t \in[0,\infty],
 \label{defnJ} \\J_{t}f(\phi_0,\phi_1) & \triangleq \inf_{s \in [t,\infty]} J
f(s,\phi_0,\phi_1),  \quad  t \in [0,\infty].
\end{align}

Recall that under $\mathbb{P}_0$, $\sigma_1$ (the first time an
observation arrives) has the exponential distribution with rate
$\mu$. Using Fubini's theorem we can write (\ref{defnJ}) as
\begin{equation}\label{eq:fubini}
Jf(t,\phi_0,\phi_1)=\int_0^{t}e^{-(\lambda+\mu)s}\left(g+\mu \cdot
S f\right)(x(s,\phi_0),y(s,\phi_1))ds, \quad t \in[0,\infty],
\end{equation}
in which $x$ and $y$ are the functions defined in
(\ref{eq:explicit-x-and-y}) and $S$ is the linear operator
\begin{equation}\label{eq:S}
S f(\phi_0,\phi_1)=
\int_{\mathbb{R}^d}f\left(\left(1+\frac{1}{\mu}\right)r(y)\phi_0
,\left(1+\frac{2}{\mu}\right)r(y)\phi_1\right) \beta_0(dy).
\end{equation}

Below we list a few useful properties of the operator $J_0$.
\begin{lemma}\label{proofJ}
For every bounded Borel function $f:\mathbb{R}_{+}^2 \rightarrow
\mathbb{R}$, the mapping $J_0 f$ is bounded. If $f$ is a concave
function, then $J_0 f$ is also a concave function. If $f_1 \leq
f_2$ are real value bounded Borel functions, then $J_0 f_1 \leq
J_0 f_2$. That is, the operator $J_0$ preserves boundedness,
concavity and ordering.
\end{lemma}

\begin{proof}
Let us define $\|f\| \triangleq \sup_{(\phi_0,\phi_1) \in
\mathbb{R}_{+}^2}|f(\phi_0,\phi_1)|<\infty$. Since $g(\cdot) \geq
g(0,0)=\lambda/c$ and $\|S(f)\| \leq \|f\|$ we can write
(\ref{eq:fubini}) as
\[
Jf(t,\phi_0,\phi_1) \geq -\left(\frac{\lambda}{c}+\mu
\|f\|\right)\int_{0}^{\infty}e^{-(\lambda+\mu)s}ds =
-\left(\frac{\lambda}{c}+\mu \|f\|\right)\frac{1}{\lambda+\mu}.
\]
Since we also have $J_{0}f(\phi_0,\phi_1) \leq
J(0,\phi_0,\phi_1)=0$, we obtain
\begin{equation}\label{eq:norm}
-\left(\frac{\lambda}{c}+\mu \|f\|\right)\frac{1}{\lambda+\mu}
\leq J_{0}f(\phi_0,\phi_1) \leq 0,
\end{equation}
 which proves the first
assertion.

The second assertion follows since $S(f)(\cdot,\cdot)$ defined in
(\ref{eq:S}) is concave if $f$ is concave, and the functions
$\phi_0 \rightarrow x(t,\phi_0)$ and $\phi_1 \rightarrow
y(t,\phi_1)$ are linear for every $t \geq 0$. The preservation of
ordering follows immediately from (\ref{eq:fubini}).
\end{proof}
\begin{cor}\label{defnofvn}
Let us define $v_n:\mathbb{R}_+^2 \rightarrow \mathbb{R}$ by
\begin{equation}\label{definevn}
v_0=0 \quad \text{and $v_n=J_0 v_{n-1}$}.
\end{equation}
  Then, for every $n \in
\mathbb{N}$, $v_n$ is bounded and concave, and $-1/c \leq v_{n+1}
\leq v_n \leq 0$. Therefore $v=\lim_{n \rightarrow \infty} v_n$,
exists, and is bounded and concave. Both $v_n$ and $v$ are
continuous (not only in the interior of $\mathbb{R}_{+}^2$), they
are increasing in each of their arguments, and their left and
right partial derivatives are bounded on every compact subset of
$\mathbb{R}_{+}^2$.
\end{cor}

\begin{proposition}\label{thm:epsilon-optimal}
For every $n \in \mathbb{N}$, $v_n$ defined in
Corollary~\ref{defnofvn} is equal to $V_n$ of (\ref{defnofVn}).
For $\varepsilon>0$, let us denote
\begin{equation}\label{rneps}
r^{\eps}_{n}(\phi_0,\phi_1) \triangleq \inf \{t\in (0,\infty]: J
v_n(t,(\phi_0,\phi_1)) \leq J_0 v_n(\phi_0,\phi_1)+\varepsilon \}.
\end{equation}
And let us define a sequence of stopping times by $S^{\eps}_1
\triangleq r_{0}^{\eps}(\p) \wedge \sigma_1$ and
\begin{equation}\label{Sne}
S^{\eps}_{n+1} \triangleq
\begin{cases}
r^{\eps/2}_n(\p) & \text{if $\sigma_1 \geq r^{\eps/2}_{n}(\p)$}
\\ \sigma_1+S_n^{\eps/2} \circ \theta_{\sigma_1}& \text{otherwise}.
\end{cases}
\end{equation}
Here $\theta_s$ is the shift operator on $\Omega$, i.e., $X_{t}
\circ \theta_s=X_{s+t}$. Then $S^{\eps}_n$ is $\eps$ optimal,
i.e.,
\begin{equation}\label{vnSne}
\E_0^{\phi_0,\phi_1}\left[\int_0^{S^{\eps}_n}e^{-\lambda
t}g(\p_t)dt \right] \leq v_{n}(\phi_0,\phi_1)+\eps.
\end{equation}
\end{proposition}

\subsection{Optimal Stopping Time}\label{optstotimesect}

\begin{proposition}\label{conoptst}
$\tau^{*}$ defined in (\ref{optst}) the smallest optimal stopping
time for (\ref{eq:value-function}).
\end{proposition}

 We will divide the proof of
this theorem into several lemmas. The following lemma shows that
if there exists an optimal stopping time it is necessarily greater
than or equal to $\tau^*$.

\begin{lemma}\label{gthants}
\begin{equation}
V(\phi_0,\phi_1)=\inf_{\tau \geq \tau^{*}}
\E_0^{\phi_0,\phi_1}\left[\int_{0}^{\tau}e^{-\lambda s }
g(\p_s)ds\right].
\end{equation}
\end{lemma}

\begin{proposition}
\label{prop:equality-of-v-and-V} We have $v(\phi_0, \phi_1)=
V(\phi_0,\phi_1)$ for every $(\phi_0,\phi_1)\in \R^2_+$. Moreover,
$V$ is the largest nonpositive solution $U$ of the equation
$U=J_0U$.
\end{proposition}

As an immediate corollary to Propositions~\ref{proofJ} and
\ref{prop:equality-of-v-and-V} and
Propositions~\ref{prop:explicit-solution}-\ref{prop:sec5-last},
which construct bounds on the optimal stopping region, we can
state the following:
\begin{cor}
Let us define the optimal stopping regions
\begin{equation}
  \bGamma_n \triangleq \{(\phi_0,\phi_1)\in \R^2_+:
  v_n(\phi_0,\phi_1)=0\}, \qquad \C_n \triangleq \R^2_+ \backslash
  \bGamma_n, \quad n\in \N,
\end{equation}
and recall that
\begin{equation}
  \bGamma = \{(\phi_0,\phi_1)\in \R^2_+:
  v(\phi_0,\phi_1)=0\},  \quad \C \triangleq \R^2_+ \backslash \bGamma.
\end{equation}
There are decreasing, convex and continuous mappings
$\gamma_n:\R_{+} \rightarrow \R_{+}$, $n \in \N$, and
$\gamma:\R_{+} \rightarrow \R_{+}$ such that
\begin{equation}
\Gamma_n=\{(\phi_0,\phi_1)\in \R_{+}: \phi_1 \geq
\gamma_n(\phi_0)\},\, \in \mathbb{N} \quad \text{and}\quad
\Gamma=\{(\phi_0,\phi_1)\in \R_{+}: \phi_1 \geq \gamma(\phi_0)\}.
\end{equation}
The sequence $\{\gamma_n(\phi_0)\}_{n\in \N}$ is increasing and
$\gamma(\phi_0) = \lim \uparrow \gamma_n(\phi_0)$ for every
$\phi_0\in \R_+$. If there are paths $t \rightarrow
(x(t,\phi_0),y(t,\phi_1))$, $t \geq 0$, $(\phi_0,\phi_1) \in \Co$,
that exit $\mathbb{C}_0$, then there exists $\xi \in
[0,\lambda/c)$ (the value of $\xi$ depends on the parameter
values) such that
$\gamma_n(\phi_0)=\gamma(\phi_0)=\lambda/c-\phi_0$ for $\phi_0
\geq \xi$, i.e., the free boundary coincides with the boundary of
the region $\mathbb{C}_0$ defined in (\ref{defn:advantageous}). In
fact if (i) $\lambda \geq 2$, or, (ii) $\lambda \in [1,2)$ and $c
\geq 2-\lambda$, or, (iii) $\lambda \in (0,1)$ and $c \geq \max
\left(2-\lambda,1-\lambda\right)$ then $\xi=0$. If (iv) $\lambda
\in [1,2)$ and $c \in (0,2-\lambda)$, (v) $\lambda \in (0,1)$ and
 $(2-\lambda)(1-\lambda)/(3-\lambda-m) < c < \max
\left(2-\lambda,1-\lambda \right)$, then
$\xi=\lambda(-1+(2-\lambda)/c)$. If on the other hand, $\lambda
\in (0,1)$ and $ c \geq 2 (1-\lambda)/(3-m)< c \leq
(2-\lambda)(1-\lambda)/(3-\lambda-m)$, then
$\xi=(2-\lambda)/c+m-3$.
\end{cor}

\begin{lemma}\label{lem:lem delay-equation}
Let $f:\R^2_+\mapsto \R$ be a bounded function. For every
  $t\in \R_+$ and $(\phi_0,\phi_1)\in \R^2_+$,
\begin{align}
\label{lem:delay-equation} J_t f (\phi_0,\phi_1)  =
Jf(t,(\phi_0,\phi_1)) +  e^{-(\lambda+\mu)t}\, J_0 f
\big(x(t,\phi_0),y(t,\phi_1))\big).
\end{align}
\end{lemma}

\begin{remark} \normalfont
\label{rem:right-continuity-of-V-Phi_t}
  Since $V$ is bounded, and $V=J_0V$ by
  Proposition~\ref{prop:equality-of-v-and-V}, we have
\begin{align}
\label{eq:delay-equation-for-V} J_t V (\phi_0,\phi_1) =
JV(t,(\phi_0,\phi_1)) + e^{-(\lambda+\mu)t}\, V
\big(x(t,\phi_0),y(t,\phi_1))\big), \quad t\in \R_+
\end{align}
for every $(\phi_0,\phi_1)\in \R^2_+$.

Let us define the $\Fb$-stopping times
\begin{align}
\label{eq:epsilon-optimal-stopping-time}
 U_{\varepsilon} \triangleq \inf\{t\ge 0: V(\wbPhi_t)\ge
 -\varepsilon\}, \qquad \varepsilon \ge 0.
\end{align}
By Remark~\ref{rem:right-continuity-of-V-Phi_t}, we have
\begin{align}
\label{eq:the-value-at-the-hitting-time}
V\big(\wbPhi_{U_{\varepsilon}}\big) \ge -\varepsilon \quad
\text{on the
  event}\quad \left\{U_{\varepsilon} <\infty\right\}.
\end{align}
\end{remark}

\begin{proposition}
\label{prop:martingale} Let $ M_t \triangleq e^{-\lambda t}
V(\wbPhi_t) + \int^{t}_0 e^{-\lambda s} g(\wbPhi_s) ds$, $t\ge 0$.
For every $n\in \N$, $\varepsilon\ge 0$, and $(\phi_0,\phi_1)\in
\R^2_+$, we have  $\E^{\phi_0,\phi_1}_0 [M_0] =
\E^{\phi_0,\phi_1}_0 [M_{U_{\varepsilon}\land \sigma_n}]$, i.e.,
\begin{align}
\label{eq:martingale}
  V(\phi_0,\phi_1) = \E^{\phi_0,\phi_1}_0 \left[ e^{-\lambda
      (U_{\varepsilon}\land \sigma_n)}V(\wbPhi_{U_{\varepsilon}\land
      \sigma_n}) + \int^{U_{\varepsilon}\land \sigma_n}_0 e^{-\lambda
      s} g(\wbPhi_s) ds \right].
\end{align}
\end{proposition}

\noindent \textbf{Proof of Proposition~\ref{conoptst}} First we
will show that $\tau^*$ is an optimal stopping time. It is enough
to show that for every $\varepsilon \ge 0$, the stopping time
$U_{\varepsilon}$
  in (\ref{eq:epsilon-optimal-stopping-time}) is an
  $\varepsilon$-optimal stopping time for the optimal stopping problem
  (\ref{eq:value-function}), i.e.,
  \begin{align*}
     \E^{\phi_0,\phi_1}_0 \left[\int^{U_{\varepsilon}}_0 e^{-\lambda
         s} g(\wbPhi_s) ds \right] \le V(\phi_0,\phi_1) + \varepsilon,
     \quad \text{for every}\quad (\phi_0,\phi_1)\in \R^2_+.
  \end{align*}
Note that the sequence of random variables
\begin{align*}
 \int^{U_{\varepsilon}\land \sigma_{n}}_0 e^{-\lambda s} g(\wbPhi_s) ds +
 e^{-\lambda (U_{\varepsilon}\land \sigma_n)}
 V(\wbPhi_{U_{\varepsilon}\land \sigma_n}) \ge - \int^{\infty}_{0}
 e^{-\lambda s}\, \frac{\lambda}{c}\, ds -\frac{1}{c}=  - \frac{2}{c}
\end{align*}
is bounded from below. By (\ref{eq:martingale}) and Fatou's Lemma,
we have
\begin{align*}
  V(\phi_0,\phi_1) &= \liminf_{n\rightarrow \infty}
  \E^{\phi_0,\phi_1}_0 \left[ \int^{U_{\varepsilon}\land \sigma_n}_0
    e^{-\lambda s} g(\wbPhi_s) ds + e^{-\lambda (U_{\varepsilon}\land
      \sigma_n)}V(\wbPhi_{U_{\varepsilon}\land \sigma_n})\right] \\
  &\ge \E^{\phi_0,\phi_1}_0 \left[\liminf_{n\rightarrow \infty}\left(
      \int^{U_{\varepsilon}\land \sigma_n}_0 e^{-\lambda s}
      g(\wbPhi_s) ds + e^{-\lambda (U_{\varepsilon}\land
        \sigma_n)}V(\wbPhi_{U_{\varepsilon}\land \sigma_n}) \right)
  \right] \\
  &= \E^{\phi_0,\phi_1}_0 \left[\int^{U_{\varepsilon}}_0 e^{-\lambda
      s} g(\wbPhi_s) ds + 1_{\{U_{\varepsilon}<\infty\}} e^{-\lambda
      U_{\varepsilon}} V(\wbPhi_{U_{\varepsilon}})\right] \\
  &\ge \E^{\phi_0,\phi_1}_0 \left[\int^{U_{\varepsilon}}_0 e^{-\lambda
      s} g(\wbPhi_s) ds \right] - \varepsilon \;
  \E^{\phi_0,\phi_1}_0 \left[ 1_{\{U_{\varepsilon}<\infty\}} e^{-\lambda
      U_{\varepsilon}} \right] \quad \text{by
    (\ref{eq:the-value-at-the-hitting-time})} \\
  &\ge \E^{\phi_0,\phi_1}_0 \left[\int^{U_{\varepsilon}}_0 e^{-\lambda
      s} g(\wbPhi_s) ds \right] - \varepsilon
\end{align*}
for every $(\phi_0,\phi_1)\in \R^2_+$. This shows that
$U_{\varepsilon}$ is an $\eps$-optimal stopping time.

Now we will show that $\tau^*$ is the smallest optimal stopping
time. Let us define
\begin{equation}
\tilde{\tau}\triangleq
\begin{cases}
\tau, & \text{if $\tau \geq \tau^{*}$},
\\ \tau+\tau^* \circ \theta_{\tau}, & \text{if $\tau<\tau^*$}.
\end{cases}
\end{equation}
Then the stopping time $\tilde{\tau}$ satisfies
\begin{equation}\label{strngmarkv}
\begin{split}
\E_{0}^{\phi_0,\phi_1}\left[\int_0^{\tilde{\tau}}e^{-\lambda s}
g(\p_s)ds \right]
 &=\E_0^{\phi_0,\phi_1}\left[\int_0^{\tau}e^{-\lambda s}g(\p_s)ds+\int_{\tau}^{\tilde{\tau}}e^{-\lambda s} g(\p_s)ds\right]
\\&=\E_0^{\phi_0,\phi_1}\left[\int_0^{\tau}e^{-\lambda s}g(\p_s)ds+e^{-\lambda \tau} \int_0^{\tau^*\circ \theta_{\tau}}
e^{-\lambda s} g(\p_{s+\tau})ds\right]
\\ &=\E_0^{\phi_0,\phi_1}\left[\int_0^{\tau}e^{-\lambda s}g(\p_s)ds+  e^{-\lambda \tau}V(\p_{\tau})\right]
\\ &\leq
\E_0^{\phi_0,\phi_1}\left[\int_0^{\tau}e^{-\lambda
s}g(\p_s)ds\right].
\end{split}
\end{equation}
Here the third equality follows from the strong Markov property of
the process $\p$. Now the proof immediately follows.\hfill
$\square$.

\subsection{Structure of the Optimal Stopping Times}
 Finally, let us describe here the structure of the
optimal stopping times. For this purpose we will need the
following lemma.
\begin{lemma}
\label{cor:r_n-as-an-exit-time} Let
\begin{align}
\label{eq:r_n} r_n (\phi_0,\phi_1) = \inf\left\{s\in(0,\infty]:
Jv_n
  \big(s,(\phi_0,\phi_1)\big)= J_0 v_n (\phi_0,\phi_1)\right\}
\end{align}
be the same as $r^{\varepsilon}_n(\phi_0,\phi_1)$ in
Proposition~\ref{thm:epsilon-optimal} with $\varepsilon = 0$. Then
\begin{align}
\label{eq:r_n-as-an-exit-time} r_n (\phi_0,\phi_1)=
\inf\left\{t>0:
  v_{n+1}\big(x(t,\phi_0),y(t,\phi_1)\big)=0 \right\}
\qquad (\inf \emptyset \equiv \infty).
\end{align}
\end{lemma}

\begin{proof} Let us fix $(\phi_0,\phi_1)\in \R^2_+$, and denote
  $r_n(\phi_0,\phi_1)$ by $r_n$. We have
  $Jv_n(r_n,(\phi_0,\phi_1))= J_0 v_n
  (\phi_0,\phi_1)=J_{r_n}v_n(\phi_0,\phi_1)$.

  Suppose first that $r_n <\infty$. Since $J_0v_{n}= v_{n+1}$, taking
  $t=r_n$ and $w=v_n$ in (\ref{lem:delay-equation}) gives
  \begin{align*}
    Jv_n(r_n,(\phi_0,\phi_1))
    = J_{r_n} v_n (\phi_0,\phi_1) =
    Jv_n(r_n,(\phi_0,\phi_1)) + e^{-(\lambda+\mu)r_n}
    v_{n+1}(x(r_n,\phi_0),y(r_n,\phi_1)).
  \end{align*}
  Therefore, $v_{n+1}(x(r_n,\phi_0),y(r_n,\phi_1))=0$.

  If $0< t < r_n$, then $Jv_n (t,(\phi_0,\phi_1))> J_0 v_n
  (\phi_0,\phi_1) = J_{r_n}v_n(\phi_0,\phi_1) =
  J_{t}v_n(\phi_0,\phi_1)$ since $u \mapsto J_u v_n (\phi_0,\phi_1)$
  is nondecreasing.  Taking $t\in (0,r_n)$ and $w=v_n$ in
  (\ref{lem:delay-equation}) imply
  \begin{align*}
    J_{0} v_n (\phi_0,\phi_1) = J_{t} v_n (\phi_0,\phi_1)=
    Jv_n(t,(\phi_0,\phi_1)) + e^{-(\lambda+\mu)t}
    v_{n+1}(x(t,\phi_0),y(t,\phi_1)).
  \end{align*}
  Therefore, $v_{n+1}(x(t,\phi_0),y(t,\phi_1))<0$ for every $t\in
  (0,r_n)$, and (\ref{eq:r_n-as-an-exit-time}) follows.

  Suppose now that $r_n = \infty$. Then we have
  $v_{n+1}(x(t,\phi_0),y(t,\phi_1))<0$ for every $t\in (0,\infty)$ by
  the same argument in the last paragraph above.  Hence, $\{t>0:
  v_{n+1}(x(t,\phi_0),y(t,\phi_1))=0\} = \emptyset$, and
  (\ref{eq:r_n-as-an-exit-time}) still holds.
\end{proof}

By Proposition~\ref{conoptst}, the set $\bGamma$ is the
\emph{optimal stopping region} for the optimal stopping problem
(\ref{eq:value-function}). Namely, stopping at the first hitting
time $U_0 = \inf\{t\in \R_+: \wbPhi_t \in \bGamma\}$ of the
process $\wbPhi=(\wPhiz,\wPhio)$ to the set $\bGamma$ is optimal
for (\ref{eq:value-function}).

Similarly, we shall call each set $\bGamma_n$, $n\in \N$ a
\emph{stopping region} for the family of the optimal stopping
problems in (\ref{defnofVn}). However, unlike the case above, we
need the first $n$ stopping regions, $\bGamma_1,\ldots,\bGamma_n$,
in order to describe an optimal stopping time for the optimal
stopping problem in (\ref{defnofVn}) (the optimal stopping times
are not hitting times of a certain set). Using
Corollary~\ref{cor:r_n-as-an-exit-time}, the optimal stopping time
$S_n \equiv S^0_n$ in Proposition~\ref{thm:epsilon-optimal} for
$V_n$ of (\ref{defnofVn}) may be described as follows: Stop if the
process $\wbPhi$ hits $\bGamma_n$ before $X$ jumps. If $X$ jumps
before $\wbPhi$ reaches $\bGamma_n$, then wait, and stop if
$\wbPhi$ hits $\bGamma_{n-1}$ before the next jump of $X$, and so
on. If the rule is not met before $(n-1)$st jump of $X$, then stop
at the earliest of the hitting time of $\bGamma_1$ and the next
jump time of $X$.

\section{Conclusion}\label{sec:conclusion}

We have solved a change detection problem for a compound Poisson
process in which the intensity and the jump size change at the
same time but the intensity changes to a random variable with a
known distribution. This problem becomes an optimal stopping
problem for a Markovian sufficient statistic. We have analyzed a
special case of this problem, in which the rate of the arrivals
moves up to one of two possible values, and the Markovian
sufficient statistic is two-dimensional, in more detail. We have
shown that the intuition that a decision would sound the alarm
only at the times when it observes an arrival does not in general
hold, see Remark~\ref{rem:only-at-jump-times}. This intuition
becomes relevant only when the disorder intensity and delay
penalty are small. Performing a sample path analysis we have been
able to find the optimal stopping time exactly for most of the
range of parameters, and tight upper and lower bounds for the rest
of the parameter range. This work has applications in insurance
risk, in which the subject Poisson  process can be viewed as the
claim arrivals process for an insurance company.


\end{document}